\documentclass[12pt]{article}
\usepackage{amsmath,amsthm,amsfonts,amssymb,latexsym}

\title{Fr\'echet-Urysohn fans in free topological groups}
\author{Taras Banakh, Du\v{s}an Repov\v{s}, Lyubomyr Zdomskyy}

\newtheorem{theorem}{Theorem}
\newtheorem{claim}[theorem]{Claim}
\newtheorem{lemma}[theorem]{Lemma}

\newtheorem{corollary}{Corollary}

\newtheorem{proposition}[theorem]{Proposition}
\newtheorem{definition}[theorem]{Definition}

\newtheorem{question}[theorem]{Question}
\newtheorem{problem}[theorem]{Problem}

\newcommand{\IZ}{\mathbb Z}
\newcommand{\IN}{{\omega}}
\newcommand{\IR}{{\mathbb{R}}}
\newcommand{\IQ}{\mathbb Q}
\newcommand{\w}{\omega}

\newcommand{\U}{\mathcal U}

\newcommand{\Auth}{\mathrm{Auth}_+}
\newcommand{\vid}{\hat{\ }}
\newcommand{\skin}{\{0,1\}^{<\omega}}

\newcommand{\td}{\tilde}
\newcommand{\I}{\mathcal I}
\newcommand{\vp}{\varphi}

\newcommand{\A}{\mathcal A}
\newcommand{\J}{\mathcal J}
\newcommand{\diam}{\mathrm{diam}}
\newcommand{\e}{\varepsilon}
\newcommand{\lef}{\mathit{left}}
\newcommand{\rig}{\mathit{right}}

\newcommand{\dd}{\mathfrak d}

\newcommand{\Top}{\mathcal{T}op}
\newcommand{\Ob}{\operatorname{Ob}}
\newcommand{\Mor}{\operatorname{Mor}}
\newcommand{\Id}{\operatorname{Id}}
\newcommand{\pr}{\operatorname{pr}}

\begin{document}

\maketitle
\begin{abstract}
   In this paper we answer the question of T.~Banakh and M.~Zarichnyi 
constructing a copy of the Fr\'echet-Urysohn fan $S_\w$ in a topological group $G$
admitting a functorial embedding  $[0,1]\subset G$. The latter means that each autohomeomorphism of $[0,1]$ extends to a continuous homomorphism of $G$.
This implies that many natural  free topological group constructions
(e.g. the constructions of the Markov free topological group, free abelian topological group, free totally bounded group, free compact group)
applied to a Tychonov space $X$ containing a topological copy of the space $\IQ$
of rationals give topological groups containing $S_\w$.
\end{abstract}
\medskip
\large

\centerline{\textbf{Introduction}}
\normalsize
\medskip

D.~Shakhmatov noticed\footnotetext{The  authors were supported
by the Slovenian Research Agency grants P1-0292-0101-04 and BI-UA/04-06-007.


\emph{2000 MSC.} 18B30, 22A05, 54C20, 54E35, 54F15.\footnotesize}
in \cite{Sh}  that the classical 
Lefschetz-N\"obeling-Pontryagin Theorem on embeddings of $n$-dimensional compacta into $\IR^{2n+1}$ 
 has no categorical counterpart: one cannot embed every finite-dimensional compact 
space $X$ into a finite-dimensional topological group $FX$ so that each 
continuous map $f:X\to Y$ extends to a continuous group homomorphism $Ff:FX\to 
FY$. The proof of this fact exploited Kulesza's example of a pathological 1-
dimensional (non-metrizable) compact space that cannot be embedded into a 
finite-dimensional topological group \cite{Kul}. However, it was discovered in 
\cite{BZa} that the problem lies already at the level of the unit interval 
$[0,1]$:  it admits no functorial embedding into any metrizable or 
finite-dimensional group. So, each topological group containing a functorially embedded 
interval is non-metrizable and thus has uncountable character. 

In light of this let us remark that the Markov free topological group $F_MI$ 
over the interval $I=[0,1]$ has character $\chi(F_MI)=\mathfrak d$ (see 
\cite{NT1}, \cite{NT4}), where $\mathfrak d$ is the well-known  uncountable 
small cardinal equal to the cofinality of the poset $(\IN^\w,\le)$. This 
cardinal is equal to the cardinality of continuum $\mathfrak c$ under Martin's 
Axiom, but can also be strictly smaller than $\mathfrak c$ in some models of 
ZFC (see \cite{vD}, \cite{Va}).

In this paper we show that the inequality $\chi(FI)\ge\mathfrak d$ holds for 
many other free topological group constructions. First we give precise 
definitions.

Let $\mathcal T$ be a subcategory of the category $\mathcal Top$ of topological spaces and 
their continuous maps. By {\em a functor of a free topological group} on 
$\mathcal T$ we understand a pair $(F,i)$ consisting of a covariant functor 
$F:\mathcal T\to\mathcal G$ from $\mathcal T$ into the category $\mathcal G$ of 
topological groups and their continuous homomorphisms, and a natural 
transformation $i=\{i_X\}:\Id\to F$ of the identity functor $\Id:\mathcal 
T\to\mathcal T$ into the functor $F$ whose components $i_X:X\to FX$ are 
topological embeddings for all spaces $X\in\mathcal T$. The naturality of $i$ 
means that for any morphism $f:X\to Y$ in $\mathcal T$ we have the following commutative 
diagram:
\begin{center}
\begin{tabular}{ccc}
 \ \  $X$       & $\stackrel{i_X}{\longrightarrow}$ &   $FX$ \ \     \\
$\ ^{f}$\Large\hskip0pt$\downarrow$\normalsize    &               & 
\Large$\downarrow$\normalsize\hskip0pt$\:^{Ff}$  \\              
\ \  $Y$  & $\longrightarrow$              &   $FY$ \ \  \\         
  &   $\;^{i_Y}$            &       
\end{tabular}
\end{center}
Therefore a functor of a free topological group $(F,i)$ to any topological 
space $X\in\Ob(\mathcal T)$ assigns a topological group $FX$ containing a 
topological copy $i_X(X)\subset FX$ of $X$ so that each morphism $f:X\to Y$ to 
another object of $\mathcal T$ extends in a canonical way to a group 
homomorphism $Ff:FX\to FY$. 
A functor $(F,i)$ is said to be \emph{minimal}, if 
for every space $X\in\mathrm{Ob}(\mathcal T)$ the group
$FX$ is algebraically
generated by $i_X(X)$.
The functor of a free compact topological group is a natural example of a non-minimal functor of a free topological group, see \cite{HM1, HM4}. 
In the case when $\mathcal T$ has only one object $X$
and $\Mor(X,X)$ coincides with the set of all autohomeomorphisms of $X$, the
embedding $i_X:X\to FX$ is simply called \cite{BZa} a \emph{functorial embedding}
of $X$ into the group $G=FX$. 
It was proven in \cite{BZa} that if there exists a functorial embedding of
the  interval $I=[0,1]$ into a topological group $G$, then $G$ is non-metrizable and 
infinite-dimensional. In this paper we shall make this result more precise by showing  that 
such a group $G$ contains a topological copy 
of the   quotient space $IS_\w=[0,1]\times\IN /\{0\}\times\IN$, 
called the
\emph{sequential hedgehog} by analogy with the sequential (or alternatively Fr\'echet-Urysohn)
fan $S_\w=S_0\times\w/\{0\}\times\w$, where $S_0=\{0\}\cup\{1/n:n>0\}$ is the convergent sequence.
This answers a question stated in \cite{BZa}, and
 implies that $\chi(G)\ge\mathfrak d$.

\begin{theorem} \label{main-additional}
If there exists a functorial embedding of the closed interval $[0,1]$ into a 
 topological group $G$, then  
$G$ contains a copy of the sequential hedgehog $IS_\w$.
\end{theorem}

In particular, topological groups fulfilling the requirements of Theorem~\ref{main-additional}  do not have the property $\alpha_4$ introduced in \cite{Ar2}.
A subcategory $\mathcal T$ of the category $\Top$  of topological 
spaces is said to be {\em full}\/ if any continuous map between objects of $\mathcal T$ is a 
morphism in $\mathcal T$. 

\begin{theorem}\label{cor} Let $(F,i)$ be a functor of a free topological 
group on a full subcategory $\mathcal T$ of $\Top$, containing the interval 
$I=[0,1]$ as an object. For every space $X\in\Ob(\mathcal T)$ containing a copy 
of $\IQ$ (resp. $I$), the topological group $FX$ contains a copy of the Fr\'echet-Urysohn 
fan $S_\w$ (resp. sequential hedgehog $IS_\w$) and hence has character $\chi(FX)\ge\dd$.
\end{theorem}

It is interesting to remark that Theorem~\ref{cor} is not true for the  
category $\mathcal T=\mathcal C(0)$ of zero-dimensional compacta and their continuous maps. 
For example, the functor $F$ which assigns to each zero-dimensional compact 
space $X$ the compact group $\IZ_2^{\Mor(X,\IZ_2)}$ containing a diagonally 
embedded copy of $X$, is a functor of a free topological group (see \cite{Sh}). For any
 compact metrizable space $X$ the 
group $\IZ_2^{\Mor(X,\IZ_2)}$ is metrizable and hence contains no copy of the 
Fr\'echet-Urysohn fan. This shows that the inclusion $I\in\Ob(\mathcal T)$ 
in  Theorem~\ref{cor} is essential.

\begin{problem} Let $(F,i)$ be a functor of a free topological group on a full 
subcategory $\mathcal I$ of $\Top$. Assume that $I\in\Ob(\mathcal T)$ and $X$ is 
an object of $\mathcal T$. Does $FX$ contain a copy of $S_\w$ if $X$ contains a 
non-trivial convergent sequence? {\rm (This is true for the  functor of the 
Markov free topological group).}
\end{problem}
We recall that a topological space $X$ is said to be \emph{scattered}, if every non-empty subspace $Y$ of $X$
has an isolated point.
Combining Theorem~\ref{cor} with the main result of \cite{Ba}, we shall derive the following
 
\begin{corollary} \label{coraa7}
Let $(F,i)$ be a minimal functor of a free topological group 
on a full category $\mathcal T$ of topological spaces such that $I\in\Ob(\mathcal T)$. 
Suppose that $X\in\Ob(\mathcal T)$ is a metrizable separable space that
has a compactification $\bar{X}\in\Ob(\mathcal T)$.
If $FX$ is a $k$-space, 
then $X$ is locally compact or scattered.
\end{corollary}

This corollary can be compared with a result \cite{AOP} of 
Arkhangel'ski\u\i, Okunev, and Pestov who proved that for a metrizable space $X$ the Markov
 free topological group  $F_MX$ is a $k$-space if and only if $X$ is either
discrete or  locally compact and separable.
\medskip

\large
\centerline{\textbf{Proof of Theorem~\ref{main-additional}}}
 \normalsize
\medskip

First  we shall describe a copy of the sequential
hedgehog $IS_\w$ in a topological group
$G$ admitting  a functorial embedding of $[0,1]$. Here the crucial role 
belongs  to special trees consisting of closed intervals and ordered
by the inverse inclusion relation, which will be called 
\emph{usual Cantor schemas} throughout the paper.

We shall use the following notations:
$\{0,1\}^{<n}=\bigcup_{k<n}\{0,1\}^k$, $\{0,1\}^{\leq n}=\{0,1\}^{<n+1}$,
 $\{0,1\}^{<\w}=\bigcup_{n\in\w}\{0,1\}^n$. In what follows 
we identify the unique element of $\{0,1\}^{0}$ with the empty sequence 
$\emptyset$.   
For  finite sequences $s=(s_0,\dots,s_n)$ and 
$t=(t_0,\ldots,t_m)$
we denote by $s\vid t$  the
\emph{concatenation}  of $s$ and $t$, i.e., the finite sequence 
$(s_0,\ldots,s_n,t_0,\ldots,t_m)$. The sequence $(0,\ldots, 0)$ of $n$ zeros
will be denoted by $0^n$. In the same way we define the sequence $1^n$.
The length $l(s)$ of a sequence $s\in\skin$ is, by definition, the number 
$n\in\w$
such that $s\in\{0,1\}^n$.

\begin{definition} A family $\J=\{J_s:s\in\skin\}$ of subsets of $\IR$
 is called a \emph{usual Cantor scheme}, if it satisfies the following conditions:
\begin{itemize}
\item[$(i)$] $J_s$ is a closed subinterval of $\IR$ for all $s\in\skin$;
\item[$(ii)$] $\min J_{s}=\min J_{s\vid 0}<\max J_{s\vid 0} < \min J_{s\vid 
1}<\max J_{s\vid 1}=\max J_s$  for all $s\in\skin$;
\item[$(iii)$] the sequence $(\diam(J_{(s_0,\ldots,s_n)}))_{n\in\w}$ converges to 
$0$
for any $(s_i)_{i\in\w}\in\{0,1\}^\w$.
\end{itemize} 
\end{definition}
For a usual Cantor scheme $\J=\{J_s:s\in\skin\}$ we make the following 
notations:
$\J_n=\{J_s:s\in\{0,1\}^n\}$, $J^{md}_s=[\max J_{s\vid 0},\min J_{s\vid 1}]$
(here \emph{md} comes from the word ``\emph{middle}''). A usual Cantor scheme 
$\J$ is called \emph{symmetric}, if the  middle points of $J_s$ and $J^{md}_s$ coincide 
for all $s\in\skin$. In what follows all usual Cantor schemas are assumed to be symmetric. 

For example, there is a canonical usual Cantor 
scheme $\I=\{I_s:s\in\skin\}$ appearing in the process of construction of
the standard Cantor set $ K\subset [0,1] $ consisting of
those numbers that have only  0's and 2's in their ternary expansion. Recall, 
that in order to obtain $K$ we exclude from  $[0,1]$ open intervals step by 
step:
 $ (\frac{1}{3},\frac{2}{3}) $ at the first step,
$(\frac{1}{9},\frac{2}{9})$ and $(\frac{7}{9},\frac{8}{9})$ at the second step, 
and so on.
Thus $I_{(0)}=[0,\frac{1}{3}]$, $I_{(1)}=[\frac{2}{3},1]$, 
$I_{(0,0)}=[0,\frac{1}{9}]$, $\ldots$

Let $ I\subset G $ be an embedding. For every $ X =\{x_{1},\ldots,x_{n}\}\subset 
I $, $ x_{1}<\cdots<x_{n} $, we set
 $\pi_-(X)=x_{1}^{-1} x_{2}\cdots x_{n}^{(-1)^{n}},$
$\pi_+(X)=x_{1} x_{2}^{-1}\cdots x_{n}^{(-1)^{n+1}},$ and
 $\pi_-(\emptyset)=\pi_+(\emptyset) =e$, where $e$ is the neutral element of 
$G$\footnote{In the rest of the paper the neutral element of a topological group $G$
will be denoted by $e_G$ or simply  by $e$ when $G$ is clear from the context.}.
For a family $\A$ of intervals we shall  denote by $\partial \A$ the set 
$\cup_{I\in\A}\partial I$,
where $\partial I$ is the set of  end-points of $I$.
We shall also write $\diam([a,b])$ for $|a-b|$. 

Given a usual Cantor scheme $\J$, we
define for every $s\in\skin$   two maps $\lef_{s,\J}, \rig_{s,\J}:[0,1]\to J_s$
such that $\lef_{s,\J}(0)=\min J_s$, $\rig_{s,\J}(0)=\max J_s$,
$$  \lef_{s,\J}(1/3^n)=\max J_{s\vid 0^{n+1}} \mbox{ and } \rig_{s,\J}(1/3^n)=\min J_{s\vid 1^{n+1}}    $$ 
for all $n\in\w$, and $\lef_{s,\J},\rig_{s,\J}$ are linear on every
interval $[1/3^{n+1}, 1/3^n]$. It is clear that $\lef_{s,\J},\rig_{s,\J}$
are continuous maps with the property 
\begin{equation} \label{left-right}
\lim_{\xi\to 0}\lef_{s,\J}(\xi)=\min J_s \mbox{ and }
\lim_{\xi\to 0}\rig_{s,\J}(\xi)=\max J_s.
\end{equation}
Whenever $\J$ is clear from the context, we shall simply write $\lef_s$ 
and $\rig_s$ in place of $\lef_{s,\J}$ and $\rig_{s,\J}$.
For every $n\in\w$ and $\xi\in (0,1]$ define the family $\J_{n,\xi}$ 
consisting of $2^{n+1}$
closed intervals as follows:
$$  \J_{n,\xi}=\big\{\lef_s([0,\xi]), \rig_s([0,\xi]) :s\in\{0,1\}^n\big\}  $$
and define the map $z_{n,\J}: [0,1]\to G$ letting $z_{n,\J}(0)=e$ and
$z_{n,\J}(\xi)=\pi_-(\partial\J_{n,\xi})$ for all $\xi\in (0,1]$.
The continuity of maps $\lef_s,\rig_s$ implies the continuity of $z_{n,\J}:[0,1]\to G$
for every $n\in\w$. Moreover, it easily follows from equation (\ref{left-right})
that $\lim_{\xi\to 0} z_{n,\J}(\xi)=e$, where $e$ is the neutral element of $G$.
Again, we shall write $z_n(\xi)$ instead of $z_{n,\J}(\xi)$ in case when $\J$ is clear
from the context.

Theorem~\ref{main-additional} is a direct consequence 
of the following technical statement,
which gives a description of a copy of the sequential hedgehog in a topological group
admitting a functorial embedding of $[0,1]$.
 
\begin{proposition} \label{descr_hedg}
If $\J$ is a usual Cantor scheme and
 $J_\emptyset\subset G$ is a functorial embedding of the closed interval $J_{\emptyset}$ into a topological group $G$, then
for every $n\in\w$ the map $z_{n}:[0,1]\to G$
is an embedding, and there exists a sequence $(d_n)_{n\in\w}\in (0,1]^\w$
such that $\bigcup_{n\in\w}z_n([0,d_n])\subset G$ is a topological copy of
the sequential hedgehog $IS_\w$.
\end{proposition}

The proof of Proposition~\ref{descr_hedg} will be split into
a sequence of lemmas. The following is the most technically difficult
among them. 

 \begin{lemma} \label{main-moroka}
If $\J$ is a usual Cantor scheme and
$J_\emptyset\subset G$ is a functorial embedding, then
$ e\not\in\overline{\bigcup_{n\in\w}z_{n,\J}([d_n,1])}$
for every   sequence $(d_n)_{n\in\w}\in (0,1]^\w$.
\end{lemma}

The proof of Lemma~\ref{main-moroka} will be also split into a sequence
of more simple lemmas.
The first of them is straightforward, and  its proof is left to the reader.
(We recall that $\I=\{I_s:s\in\{0,1\}^{<\w}\}$ is the canonical Cantor scheme.)

\begin{lemma} \label{pr2006}
Let $\J$ be a usual Cantor scheme. Then 
there exists an increasing homeomorphism  $h:J_\emptyset\to [0,1]$ such that
$h(J_s)=I_s$ for all $s\in\skin$.
\end{lemma}

The subsequent lemma easily   follows from the  compactness of $ [0,1] $.

\begin{lemma} \label{aaa} 
Let $[0,1]\subset G$ be an embedding and $U$ be an open neighborhood 
of the neutral element $e\in G$. Then
there exists $\varepsilon (U)>0 $ and an open neighborhood $V$ of $e$
 such that for every
$ x_{1}, x_{2}\in [0,1] $ with $|x_1-x_2|<\varepsilon(U)$ the following inclusion 
holds:
\[ \{x_{1}^{-1}V x_{2}, x_1 V x_2^{-1}, x_2^{-1}V x_1, x_2 V x_1^{-1}\} \subset U \]
\end{lemma}

For every $ n>0$ and $ f:X\rightarrow X$ we shall denote 
 by $f^n$ the composition $\underbrace{f\circ\cdots\circ f}_{n}$. It 
will be also convenient for us to denote by $f^0$
the identity map on $X$. 
Given any strictly increasing   function $ \vp\in\w^\w $ such that $\vp(0)>0$
 and a usual Cantor scheme $\J$, define a usual 
 Cantor scheme $\td{\J}=\{\td{J}_s:s\in \{0,1\}^{<\w}\}$
letting $\td{J}_{\emptyset}=J_\emptyset$ and
 $$ \td{J}_{(s_0,\ldots,s_{n})}=J_{(s_0^{\vp^2(0)}\vid s_1^{\vp^4(0)-\vp^2(0)} 
\vid \cdots \vid s_n^{\vp^{2n+2}(0)-\vp^{2n}(0)} )}.$$ 
Thus $\td{\J}_m\subset\J_{\vp^{2m}(0)}$ for all $m\in \w$.

\begin{claim} \label{struc}
Let $p\in \w$ and $\vp^{2p}(0)\leq n<\vp^{2p+1}(0)$. Then 
$\partial \J_{n,\xi}\subset\{0,1\}\cup\bigcup_{s\in\{0,1\}^{\leq p}}\td{J}^{md}_s$
for all $\xi> 3^{n-\vp(n)}$.
\end{claim}
\begin{proof}
Let us fix arbitrary $a\in\partial\J_{n,\xi}$ and $\xi>3^{n-\varphi(n)}$. 
Assume, contrary to our claim, that $a$ lies in the interior  of
$\td{J}_s$  for
some $s\in\{0,1\}^{p+1}$. Concerning $a$, two cases are possible.

1. $a\in\{\lef_t(0),\rig_t(0)\}$ for some $t\in\{0,1\}^n$. Then $a\in\partial J_t$.
Since $n<\vp^{2p+2}(0)$ and $\td{J}_s\in \J_{\vp^{2p+2}(0)}$, the inclusion
$a\in\mathrm{Int}(\td{J}_s)$ is impossible.

2. There exists $t\in\{0,1\}^n$ such that $a\in\{\lef_t(\xi),\rig_t(\xi)\}$. 
Without loss of generality, $a=\lef_t(\xi)$. Set $b=\min J_t=\lef_t(0)$.
As it was already proven in the case~1, 
$b\in \{0,1\}\cup\bigcup_{s\in\{0,1\}^{\leq p}}\td{J}^{md}_s$. Again, two subcases are possible:

$(i)$\   $b$ is an end-point of $\td{J}^{md}_{s_0}$ for some $s_0\in\{0,1\}^{\leq p}$,
or $b=0$. Then there exist $s_1\in\{0,1\}^p$ and $i\in\{0,1\}$ such that
 $b\in\partial\td{J}_{s_1\vid i}$.
Since $b=\min J_t$, we conclude that $i=0$ and $b=\min\td{J}_{s_1\vid 0}$.
Let $r\in\{0,1\}^{\vp^{2p}(0)}$ be such that $\td{J}_{s_1}=J_r$.
The inequality $n\geq\vp^{2p}(0)$ combined with $J_t\in\J_n$, $J_r\in\J_{\vp^{2p}(0)}$,
  and $\min J_r =\min J_t=b$,
implies that 
$$a=\lef_t(\xi)\leq\lef_r(\xi)\leq\lef_r(1)=\max J_{r\vid 0}<\min J_{r\vid 1}\leq
\min \td{J}_{s_1\vid 1}.$$
In addition,
$$a=\lef_t(\xi)>\lef_t(3^{n-\vp(n)})=\max J_{t\vid 0^{\vp(n)-n+1}}\geq
\max J_{t\vid 0^{\vp^{2p+2}(0)-n}}= \max \td{J}_{s_1\vid 0}. $$
(The last equality immediately follows from  $b =\min J_{t\vid 0^{\vp^{2p+2}(0)-n}}=
\min \td{J}_{s_1\vid 0} $ and $J_{t\vid 0^{\vp^{2p+2}(0)-n}},$ 
$\td{J}_{s_1\vid 0}\in\J_{\vp^{2p+2}(0)}$.) Therefore
$a\in (\max\td{J}_{s_1\vid 0},\min\td{J}_{s_1\vid 1})=\mathrm{Int}(\td{J}^{md}_{s_1})$,
a contradiction.

$(ii)$\   $b\in\mathrm{Int}\td{J}^{md}_{s_2}$ for some $s_2\in\{0,1\}^{\leq p}$.
Again, let
  $r\in\{0,1\}^{\vp^{2p}(0)}$ be such that $\min J_r=\min J_t = b$.
Since $\min J_r\in\mathrm{Int}(\td{J}^{md}_{s_2})$, we conclude that
$J_r\subset\mathrm{Int}(\td{J}^{md}_{s_2})$.
The inequality $n\geq\vp^{2p}(0)$ combined with  $\min J_r =\min J_t=b$
implies that $J_t\subset J_r\subset \mathrm{Int}(\td{J}^{md}_{s_2})$.
Therefore $a\in \mathrm{Int}(\td{J}^{md}_{s_2})$ being an element of $J_t$,
which contradicts the equation $a=\mathit{left}_t(\xi)$.
\end{proof}

In the sequel the notation $\Auth([a,b])$ stands for the family of all 
increasing
autohomeomorphisms of the interval $[a,b]$. 

\begin{lemma} \label{tech1}
Let $\vp$ and $\td{\J}$ be  as above. Then 
for every indexed family $\{U_s:s\in\{0,1\}^{<\w}\}$ of open neighborhoods of 
the identity
$e$ of $G$ there exists a sequence
 $(h_n)_{n\in\w}\subset\Auth([0,1])$ such that 
\begin{itemize}
\item[$(1)$] $h_n(\td{\J})$ is a symmetric usual Cantor scheme; 
\item[$(2)$]$h_{n+1}|\bigcup_{s\in\{0,1\}^{\leq n}}\td{J}^{md}_s=
h_n|\bigcup_{s\in\{0,1\}^{\leq n}}\td{J}^{md}_s$ for all $n\in\w$; and
\item[$(3)$] $\{ \pi_{-}[h_n(\partial \J_{m,\xi} \cap\td{J}^{md}_s)], 
\pi_{+}[h_n(\partial \J_{m,\xi} \cap\td{J}^{md}_s)]\}  \subset U_s$
for all $s\in\{0,1\}^{\leq n}$, $m\geq \vp^{2n}(0)$, and $\xi>0$.
\end{itemize}
\end{lemma}
\begin{proof}
The following claim
is the main building block of our proof.
In order to shorten its proof, some explanations similar to those made in the proof
of Claim~\ref{struc} are omitted.
\begin{claim} \label{tech0}
Let $\J$, $G$, $e$ be  as in Lemma~\ref{tech1}.
 Let also $s\in\{0,1\}^{n_0}$, $m>0$,
 $c=\max J_{s\vid 0^m}$, $d=\min J_{s\vid 1^m}$, and let 
$U$ be an open neighborhood of $e$ in $G$.
Then there exists $h\in\Auth(J_s)$
such that 
\begin{itemize}
\item[$(i)$] $h|J_{s\vid 0^m}$ and $h|J_{s\vid 1^m}$ are linear;
\item[$(ii)$] $\diam(h(J_{s\vid 0^m}))=\diam(h(J_{s\vid 1^m}))$; and
\item[$(iii)$] $\{\pi_-(h(\partial \J_{n,\xi}\cap [c,d])), \pi_+(h((\partial 
\J_{n,\xi}\cap [c,d]))\}\subset U$
for all $n\geq n_0$ and $\xi\in (0,1]$.
\end{itemize}
\end{claim}
\begin{proof}
Let us find some $a,b\in J_s$ such that $a<b$, the middle points of $[a,b]$ and $J_s$ coincide,
and $|a-b|<\varepsilon(U)$. 
 The latter means that there exists an open neighborhood $V$ of $e$  
such that $u^{-1}V v\cup u V v^{-1}\subset U$ for all $u,v\in [a,b]$.
The continuity of the group operation on $G$ gives us a sequence
$(V_n)_{n\in\w}$ of open neighborhoods of $e$ such that $V_n^{2^n}\subset V$.
Let  $h\in\Auth(J_s)$ be such that the following conditions are satisfied:
$h|[\min J_s,c]$ and $h|[d,\max J_s]$ are linear bijections
onto $[\min J_s,a]$ and $[b,\max J_s]$ respectively, and  
$\diam(h(J_t))<\e(V_{l(t)})$ provided  $J_t\subset (c,d)$.
The existence of $h$ follows from Lemma~\ref{pr2006}.
Given arbitrary $n\geq n_0$ and $\xi\in [0,1]$, write the family
$\{J\in\J_{n,\xi}:J\subset [c,d]\}$ in the form
$\{J_1,\ldots,J_q\}$   such that $J_i$ is situated to the left of
$J_j$ provided $ i<j$. Let us note that each $J\in\J_{n,\xi}$
is contained in some element of $\J_{n+1}$, and hence $q\leq 2^{n+1}$.
It can be easily derived from the definitions of a usual Cantor scheme
and maps $\lef_{-},\rig_{-}$ 
that  $\{c,d\}\subset\partial \J_{n,\xi}$ for all $n,\xi$ such that
$ \xi=3^{-(n_0+m-n-1)}$ 
or $n\geq n_0+m $,
 and $\{c,d\}\cap\partial\J_{n,p}=\emptyset$ otherwise.
In the first case we have
\begin{eqnarray*}
\pi_{-}[h(\partial \J_{n,\xi}\cap [c,d])] =a^{-1}
\pi_+(h(\partial J_1))\pi_+(h(\partial J_2))\cdots\pi_+(h(\partial J_q)) 
b\subset 
 a^{-1}V_{n+1}^q b\subset a^{-1}V b\subset U.  
\end{eqnarray*}  
In the second case $n<n_0+m$ and $\xi\neq 3^{-(n_0+m-n-1)}$.
If $\xi<3^{-(n_0+m-n-1)}$, then $(c,d)$ contains all elements of $\J_{n,\xi}$
whose intersection with $[c,d]$ is nonempty,  
and hence
\begin{eqnarray*}
\pi_{-}[h(\partial \J_{n,\xi}\cap [c,d])] = 
\pi_-(h(\partial J_1))\pi_-(h(\partial J_2))\cdots\pi_-(h(\partial J_q)) \subset 
 V_{n+1}^q \subset U.  
\end{eqnarray*}  
It sufficies to consider the case $\xi>3^{-(n_0+m-n-1)}$. Let $u=\lef_{s}(\xi)$
and $v=\rig_s(\xi)$. Then $c<u<v<d$ and $J_i\subset (u,v)$
for all $i\leq q$. Therefore
\begin{eqnarray*}
\pi_{-}[h(\partial \J_{n,\xi}\cap [c,d])] =
u^{-1}\pi_+(h(\partial J_1))\pi_+(h(\partial J_2))\cdots\pi_+(h(\partial J_q)) v 
\subset 
 u^{-1} V_{n+1}^q v \subset u^{-1}V v\subset U.  
\end{eqnarray*} 

Verification that $\pi_+(h(\partial \J_{n,\xi}\cap [c,d]))\subset U$ is similar,
and thus condition $(iii)$ is satisfied. 
\end{proof}

Applying Claim~\ref{tech0} for the usual Cantor scheme $\J$,
 $s=\emptyset$, $m=\vp^{2}(0)$, and $U=U_\emptyset$, we get 
 $h_0\in\Auth(\td{J}_\emptyset)$ satisfying the conditions $(i)-(iii)$  of 
Claim~\ref{tech0}.  Conditions $(i)$ and $(ii)$ obviously imply $(1)$,
condition $(iii)$ implies $(3)$, while $(2)$ is trivial.

Assuming that we have already constructed $h_k$ satisfying $(1)-(3)$
 for all $k<n$,
set $h_n|\bigcup_{s\in\{0,1\}^{\leq n-1}}\td{J}_s^{md}=h_{n-
1}|\bigcup_{s\in\{0,1\}^{\leq n-1}}\td{J}_s^{md}$.
Thus condition $(2)$ is satisfied.
In addition,  $(3)$ holds for all $s\in\{0,1\}^{\leq n-1}$.
It suffices to define $h_n$ on  $[0,1]\setminus\bigcup_{s\in\{0,1\}^{\leq n-
1}}\td{J}_s^{md}=\bigcup_{s\in\{0,1\}^n}(\td{J}_s\setminus\partial\td{J}_s)$. 
Let  us note, that for every particular $s\in\{0,1\}^n$ the construction of 
$h_n|\td{J}_s$
is similar to that of $h_0$.
Given any $s=(s_0,\ldots,s_{n-1})\in\{0,1\}^n$, set $t=s_0^{\vp^2(0)}\vid 
s_1^{\vp^4(0)-\vp^2(0)}\vid\cdots\vid s_{n-1}^{\vp^{2n}(0)-\vp^{2n-2}(0)}$ and
$m=\vp^{2n+2}(0)-\vp^{2n}(0)$. Thus $J_t=\td{J}_s$.
Applying Claim~\ref{tech0} for $\J$,
 $t\in\{0,1\}^{\vp^{2n}(0)}$, $m$, and $U_s$, we get $h_s\in\Auth(J_t)$
satisfying conditions $(i)-(iii)$ of Claim~\ref{tech0}. Set 
$h_n|\td{J}_s=h_s\circ h_{n-1}|\td{J}_s$. Again, $(i)$ and $(ii)$ imply $(1)$, 
and $(iii)$ implies $(3)$.
\end{proof}

The following simple statement is due to M.~Tkachenko (see 
\cite[Lemma~1.3]{Tk-ch}). 
\begin{lemma} \label{seq} 
Let $G$ be a topological group and $\U$ be the family of all open neighborhood 
of the neutral element of $G$. Then for every   $ U\in\mathcal{U}$
there exists a decreasing sequence $(U_n)_{n\in\omega}\subset\mathcal{U}$ such 
that 
$U_{\sigma(0)}U_{\sigma(1)}\cdots U_{\sigma(n)}\subset U$ for every bijection
$\sigma:\{0,\ldots,n\}\to\{0,\ldots,n\}$.
\end{lemma}

The proof of the following simple technical lemma is left to the reader.
\begin{lemma} \label{bbb} 
For every function $ \vp:\omega\rightarrow\w\setminus\{0\}$ there exist strictly increasing 
functions
$\psi, \theta : \omega\rightarrow\IN $ such that
 \begin{itemize} 
\item[$(i)$] $\varphi(k)<\psi(k)<\theta (k)$  for all $k\in\omega$;  and
\item[$(ii)$] $\psi^{n+1}(0)=\theta^n(0)$  for all $n\in\IN$.
\end{itemize} 
\end{lemma}
\medskip

\noindent \textbf{Proof of Lemma~\ref{main-moroka}}.
Throughout the proof we denote by $\hat{h}:G\to G$ some continuous homomorphism
extending $h\in\Auth(J_\emptyset)$. 
Let $\U$ be the family of all open neighborhoods of the neutral element $e$ of 
$G$. The continuity of the group operation of $G$ yields the existence of an element 
$ U\in\mathcal{U}$ such that $ U\cap 0^{-1} U\: 1 =\emptyset $. By Lemma~\ref{seq}
there exists a sequence $(U_{n})_{n\in\omega}\subset\mathcal{U}$
such that $U_{i_1}U_{i_2}\cdots U_{i_{2^{n+1}-1}}\subset U$
for any $n\in\w$ and $(i_1,\ldots,i_{2^{n+1}-1})$
such that $|\{j<2^{n+1}:i_j=k\}|=2^k$ for all $k\in\{0,\ldots, n\}$.

Let $\vp:\w\to\w$ be a strictly increasing map with the property 
$3^{n-\vp(n)}<d_n$ for every  $n\in\w$. 
Let us fix a sequence $(h_n)_{n\in\w}\subset\Auth([0,1])$ satisfying conditions
$(1)-(3)$ of Lemma~\ref{tech1} with $U_s$ equal to   $U_{l(s)}$ defined above, 
where $s\in\{0,1\}^{<\w}$.
Conditions $(1)$, $(2)$ imply that
$|h_n(t)-h_{n+1}(t)|\leq 2^{-n-1}$, and hence the sequence $(h_n)_{n\in\w}$
is uniformly convergent to a monotone continuous surjective function $h:[0,1]\to 
[0,1]$.
For every $s\in\{0,1\}^{\leq n}$ we have $h|\td{J}_s^{md}=h_n|\td{J}_s^{md}$ by 
$(2)$,
consequently    $h|\td{J}_s^{md}$ is not constant for every $s\in\skin$.
Since $\bigcup_{s\in\skin}\td{J}_s^{md}$ is dense in $[0,1]$, we conclude that
$h\in\Auth([0,1])$ as a monotone continuous surjection which fails to be 
constant on arbitrary open subset of $[0,1]$.  We claim that 
$$\hat{h}(\bigcup_{p\in\w}\bigcup_{\vp^{2p}(0)\leq n < \vp^{2p+1}(0)}\{z_{n,\J}(\xi): 
\xi> 3^{n-\vp(n)} \})\: \subset 0^{-1}\: U\: 1.$$
Indeed, let us fix arbitrary $p\in\w$, $\vp^{2p}(0)\leq n <\vp^{2p+1}(0)$,
and $\xi\in (3^{n-\vp(n)},1]$.
Then by Claim~\ref{struc}  we have 
$\partial \J_{n,\xi}\subset \{0,1\}\cup\bigcup_{s\in\{0,1\}^{\leq p}}\td{J}_s^{md}.$
Therefore
$$  z_{n,\J}(\xi)=\pi_-(\partial\J_{n,\xi})=\pi_-(\{0\}\cup\bigcup_{s\in\{0,1\}^{\leq 
p}}(\partial \J_{n,\xi}\cap\td{J}_s^{md})\cup\{1\}).$$
Combining the equation above  with $(2)$ and $(3)$ of Lemma~\ref{tech1},
we conclude that 
\begin{eqnarray*}
\hat{h}(z_{n,\J}(\xi))=\hat{h}[\pi_{-}( \{0\}\cup\bigcup_{s\in\{0,1\}^{\leq p}}(\partial 
\J_{n,\xi}\cap\td{J}_s^{md})\cup\{1\})]=\\
=\pi_{-}[h ( 
\{0\}\cup\bigcup_{s\in\{0,1\}^{\leq p}}(\partial 
\J_{n,\xi}\cap\td{J}_s^{md})\cup\{1\})]=\\
= \pi_{-}(\{0\}\cup h_{p}(\partial\J_{n,\xi}\cap\td{J}_{s_1}^{md})\cup 
h_{p}(\partial\J_{n,\xi}\cap\td{J}_{s_2}^{md})\cup \cdots\cup  
h_{p}(\partial\J_{n,\xi}\cap\td{J}_{s_{2^{p+1}-1}}^{md})\cup \{1\})= \\
=0^{-1}\: \pi_{\delta_1}(h_{p}(\partial\J_{n,\xi}\cap\td{J}_{s_1}^{md})) 
\pi_{\delta_2}    (h_{p} (\partial\J_{n,\xi}\cap\td{J}_{s_2}^{md})) \cdots 
\pi_{\delta_{2^{p+1}-1}}(h_{p}   (\partial\J_{n,\xi}\cap\td{J}_{s_{2^{p+1}-
1}}^{md}))\:1\subset\\
\subset 0^{-1}\:U_{l(s_1)}U_{l(s_2)}\cdots U_{l(s_{2^{p+1}-1})}\:1\subset 0^{-
1}\:U\:1,   
\end{eqnarray*}
where $\{s_1,\ldots,s_{2^{p+1}-1}\}$ is the enumeration of $\{0,1\}^{\leq p}$
such that $\td{I}^{md}_{s_i}$ is situated to the left of $\td{I}^{md}_{s_j}$ 
provided
$1\leq i<j<2^{p+1}$, and $\delta_i\in\{+,-\}$. Then 
$\bigcup_{p\in\w}\bigcup_{\vp^{2p}(0)\leq n <\vp^{2p+1}(0)}\{z_{n,\J}(\xi):1\geq \xi 
 > 3^{n-\vp(n)}\} \cap \hat{h}^{-1}(U)=\emptyset$
by our choice of $U$. 

Let $\psi,\theta:\w\to\w$ be an increasing number sequences such as in 
Lemma~\ref{bbb}, i.e. $\vp(n)\leq \psi(n)\leq\theta(n)$ for all $n\in\w$ and 
$\theta^n(0)=\psi^{n+1}(0)$ for all
$n\geq 1$. It follows from the above  that there are continuous homomorphisms 
$h_\psi,h_\theta: G\to G$ such that 
\begin{eqnarray*}
 \bigcup_{p\in\w}\bigcup_{\psi^{2p}(0)\leq n <\psi^{2p+1}(0)}\{z_{n,\J}(\xi):1\geq \xi 
> 3^{n-\psi(n)}\} \cap\  h_\psi^{-1}(U)=\emptyset\  \mbox{and} \\
   \bigcup_{p\in\w}\bigcup_{\theta^{2p}(0)\leq n <\theta^{2p+1}(0)}\{z_{n,\J}(\xi):1\geq 
\xi > 3^{n-\theta(n)}\}\ \cap\  h_\theta^{-1}(U)=\emptyset,
\end{eqnarray*}
which implies
$ \bigcup_{n\in\w}\{z_{n,\J}(\xi):1\geq \xi > 3^{n-\vp(n)}\} \cap (h_\psi^{-1}(U)\cap 
h_\theta^{-1}(U))=\emptyset,$
and hence the fact that $e\not\in \overline{\bigcup_{n\in\w}\{z_{n,\J}(\xi):1\geq 
\xi> d_n\}}$ is proven.
\hfill $\Box$
\medskip

For a subset $X$ of a topological group $G$ we shall denote by
$\langle X\rangle$ the smallest subgroup of $G$ containing $X$.
\begin{lemma} \label{injection}
Let $[0,1]\subset G$ be a functorial embedding, 
$\{x_i:0\leq i\leq n\}\cup\{y_j:0\leq j\leq m\}\subset [0,1]$,
and $y_{j_0}\not\in\{y_j:j\neq j_0\}\cup\{x_i:0\leq i\leq n\}$ for some $j_0$.
Then
$$x_0^{k_0}\cdot x_1^{k_1}\cdot\cdots\cdot x_n^{k_n}\neq y_0^{l_0}\cdot y_1^{l_1}\cdot\cdots\cdot y_m^{l_m}$$
for arbitrary integers $k_i,l_i$ such that $l_{j_0}\in\{-1,1\}$.

In particular,
for every $n\in\w$  and usual Cantor scheme $\J$ the map $z_{n,\J}:[0,1]\to G$
is an embedding, and $z_{n,\J}((0,1])\cap z_{m,\J}((0,1])=\emptyset$ for all $m\neq n$.
\end{lemma}
\begin{proof}
The second statement easily follows from the first one.
To prove the first statement, set $x=x_0^{k_0}\cdot x_1^{k_1}\cdot\cdots\cdot x_n^{k_n}$, 
$y=y_0^{l_0}\cdot y_1^{l_1}\cdot\cdots\cdot y_m^{l_m}$, and
assume to the contrary that $x=y$.
Let $h\in\Auth([0,1])$ be such that $h(u)=u$ for all $u\in \{y_j:j\neq j_0\}\cup\{x_i:0\leq i\leq n\}$ and $h(y_{j_0})\neq y_{j_0}$. Since the embedding $[0,1]\subset G$ is functorial,
there exists a continuous homomorphism $\hat{h}:G\to G$  extending $h$.
It follows from the above  that $\hat{h}(x)=x=y$ and $\hat{h}(y)\neq y$,
a contradiction.
\end{proof}

\noindent\textbf{Proof of Proposition~\ref{descr_hedg}.}
In light of Lemmas~\ref{main-moroka} and \ref{injection} we are left with the task of constructing a sequence $(d_n)_{n\in\w}$ of positive reals with the property
$z_n((0,d_n])\cap\overline{\bigcup_{k>n}z_k([0,d_k])}=\emptyset$
for all $n\in\w$. Let $A=\{\xi\in (0,1]:z_0(\xi)\in\overline{\bigcup_{k>0}z_k([0,1]))}\}$.
Since $\lim_{\xi\to 0}z_k(\xi)=e$ for all $k\in\w$,
we conclude that $A$ consists of  $\xi\in [0,1]$ such that 
$ \xi\in\overline{\bigcup_{k>0}z_k([c_k,0])} $ 
for some sequence $(c_k)_{k>0}$ of positive reals. 

We claim that $0\not\in\overline{A}$.
Indeed, assuming the converse we could find a sequence $(\xi_n)_{n\in\w}$ of elements
of $A$ converging to $0$. It follows from the above  that for every $n\in\w$
there exists a sequence $(c_{n,k})_{k>0}$  of positive reals such that $z_0(\xi_n)\in\overline{\bigcup_{k>0}z_k([c_{n,k},1])}.$
Set $c_k=\min\{c_{n,k}:n\leq k\}$. Then $e\in\overline{\bigcup_{k>0}z_k([c_k,1])},$
which contradicts Lemma~\ref{main-moroka}.
 
It follows from the above that $A\subset (d_0,1]$ for some $d_0>0$,
and consequently $z_0((0,d_0])\cap\overline{\bigcup_{k>0}z_k([0,1])}=\emptyset$.
 In the same way
for every $n>0$ we can find $d_n$ such that
$z_n((0,d_n])\cap\overline{\bigcup_{k>n}z_k([0,1])}=\emptyset$,
which completes our proof.
\hfill $\Box$
\medskip

\noindent\textbf{Remark.}
Let $\J$ be a usual Cantor scheme
such that $J_\emptyset\subset G$ is a functorial embedding.
Let also $C_n=z_{n,\J}([0,d_n])$ for a sequence $(d_n)_{n\in\w}$
fulfilling the requirements of Proposition~\ref{descr_hedg}.
Then it can be easily derived from Lemma~\ref{injection} that
the map 
$$ \prod_{i\leq n}[0,d_n] \ni(\xi_0,\ldots,\xi_n)\mapsto \prod_{i\leq n}z_{n,\J}(\xi_i) $$
is a homeomorphism, and we denote its image by $D_n$.
Thus we have an increasing sequence 
$$D_0\subset D_1\subset\cdots\subset D_n\subset\cdots,$$
where each $D_n$ is a homeomorphic copy of the $(n+1)$-dimensional cube\footnote{This gives an alternative
(but much longer) proof of the fact \cite{BZa} that each topological group is infinite-dimensional provided it admits a functorial embedding
of the closed interval.}.
Let $\nu$ be the topology of $G$ and let
$\tau$ be the strongest topology on $D=\bigcup_{n\in\w}D_n$
such that $\tau|D_n=\nu|D_n$ for all $n\in\w$.
It is easy to see that
$(D,\tau)$ contains a  homeomorphic copy of $\IR^\infty=\lim_{\to}\IR^n$, which is homeomorphic to
the Markov free topological group  over $[0,1]$ (see \cite{Za}). But we do not
know whether $\nu|D=\tau|D$. This leads us to the following question.
\begin{question}
 Let $[0,1]\subset G$ be a functorial embedding.
Does $G$ contain a topological copy of the Markov free topological group over
$[0,1]$? More precisely, does the construction  described above yield a
copy of $\IR^\infty$ in $G$? 
\end{question}

Similarly to the proof of Proposition~\ref{main} below,
the positive solution of the above question would imply that for every functor
$F$ of a free topological group and every Tychonov space $X$ containing a topological copy of
$[0,1]$, the group $FX$ contains a topological copy of the Markov free topological group
over $[0,1]$.
\hfill $\Box$

\medskip
\large

\centerline{\textbf{Theorem~\ref{cor} and its generalization}}
\normalsize
\medskip

We shall derive Theorem~\ref{cor} from the following slightly more general statement, where we specify the properties of the
category  $\mathcal T$ used in the proof of Theorem~\ref{cor}.

\begin{proposition}\label{main} Let $(F,i)$ be a functor of a free topological group 
on a category $\mathcal T$ of topological spaces such that $I\in\Ob(\mathcal T)$ 
and $\Mor(I,I)\supset\Auth([0,1])$. The group $FX$ over an object $X$ of $\mathcal T$
 contains: 
\begin{itemize}
\item[$(1)$]  a copy of the Fr\'echet-Urysohn fan $S_\w$ provided there is a 
morphism $f:X\to I$ in $\mathcal T$ whose restriction $f|Q$ onto some subspace 
$Q\subset X$ without isolated points is an embedding of $Q$ into $I$,
and $\Mor(I,I)$ contains all continuous maps $\phi:I\to I$; and
\item[$(2)$] a copy of the sequential hedgehog  $IS_\w$ provided there is a 
copy $Y\subset X$ of $I$, a surjective map $f\in\Mor(X,Y)$, and
a  homeomorphism $h\in\Mor(Y,I)$.
\end{itemize}
\end{proposition}

Proposition~\ref{main} is a consequence of Theorem~\ref{main-additional}
and the following 

\begin{lemma} \label{tiur}
Let $X$ be a Tychonov space containing a topological copy $Q$ of
the space  of rational numbers and $\J$ be a usual Cantor scheme with $J_\emptyset=[0,1]$.
Then there exists 
a homeomorphic copy $Q'\subset Q$ of $Q$,
and a continuous map $\psi:X\to J_\emptyset$ such that $\psi|Q'$ is a homeomorphism
between the spaces $Q'$ and $\partial \J=\bigcup_{s\in\skin}\partial J_s$.
\end{lemma}
\begin{proof}
Let $I=[0,1]$.
It is well-known that the diagonal 
map $\delta:X\to I^{\Mor(X,I)}$ is an embedding, where $\Mor(X,I)$ stands for 
the set of all continuous maps from $X$ to $I$. Assume that $X$ contains a 
subset $Q\subset X$ homeomorphic to the space of rational numbers. Using the 
fact that $\delta(Q)\subset I^{\Mor(X,I)}$ has a countable base, one can construct 
a countable subset $C\subset \Mor(X,I)$ such that the restriction 
$\pr|\delta(Q)$ of the projection $\pr:I^{\Mor(X,I)}\to I^C$ is a homeomorphic 
embedding.  By the standard argument (see \cite[Theorem~21.18]{Ke}), we can find a 
topological copy $Q_1\subset \pr\circ\delta(Q)$ of $\IQ$ whose closure 
$\overline{Q_1}$ in the (metrizable) cube $I^C$ is zero-dimensional,
and hence is homeomorphic to the Cantor space $\{0,1\}^\w$. 
Let 
$$e:\overline{Q_1}\to \bigcup_{(s_n)_{n\in\w}\in\{0,1\}^\w}\bigcap_{n\in\w}J_{(s_0,\ldots,s_n)}$$ 
be the homeomorphism with the property $e(Q_1)=\partial\J$
(its existence follows from \cite[Part~4,~Th.~1]{BP} or the main result of
\cite{H}),
and $\bar e:I^C\to I$ be an extension 
of $e$ to a continuous map. Then  the set 
$Q'=Q\cap (\pr\circ\delta)^{-1}(Q_1)$ is a topological copy of $\IQ$ and the map $f=\bar e\circ \pr\circ \delta:X\to I$ restricted to $Q'$ is a homeomorphism between $Q'$  
and $\partial\J$.
\end{proof}

Observe that for an arbitrary family $\{(x^n_k)_{k\in\w}:n\in\w\}$ of sequences of elements
of $(0,1]$, the subspace $\{x^n_k\times\{n\}:n,k\in\w\}\cup\{\{0\}\times\w\}$ of $IS_\w$
is homeomorphic to $S_\w$ provided $\lim_{k\to\infty}x^n_k=0$ for all $n\in\w$.
\medskip

\noindent\textbf{Proof of Proposition~\ref{main}.}
Throughout the proof we shall identify $X$ with $i_X(X)$.
We present here only the proof of the first part. The proof of the second one is analogous.

 Let $f\in\Mor(X,I)$ be such that $f|Q$ is an embedding for some
homeomorphic copy $Q\subset X$ of $\IQ$.
Lemma~\ref{tiur} imples that there exists a usual Cantor scheme $\J$,
a continuous map $g:I\to I$, and a copy $Q_1\subset f(Q)$ of $\IQ$
such that $g(Q_1)=\partial\J$ and $g|Q_1$ is an embedding. Then $h=g\circ f$ belongs to
$\Mor(X,I)$ and $h(f^{-1}(Q_1))=\partial\J$.
Set $r_{n,k}=z_{n,\J}(3^{-k})$. Applying Proposition~\ref{descr_hedg},
we can find a sequence $(k_n)_{n\in\w}$ of natural numbers
such that $R=\{e\}\cup\{r_{n,k}:n\in\w,k\geq k_n\}$ is a homeomorphic copy of
$S_\w$. By our construction of the maps $z_{n,\J}$,
for every $n,k\in\w$ we can find elements
$$ 0=u_{n,k,0}<u_{n,k,1}<\cdots< u_{n,k,2^{n+2}-1}=1 $$
of $\partial\J$ such that $r_{n,k}=u_{n,k,0}^{-1}u_{n,k,1}u_{n,k,2}^{-1}\cdots u_{n,k,2^{n+2}-1}$.
In addition, $u_{n,k,p}$ does not depend on $k$ provided $4$ divides $p$ or $p+1$,
and
\begin{equation} \label{convpre}
\lim_{k\to\infty}u_{n,k,4q+1}=u_{n,k,4q}, \  \lim_{k\to\infty}u_{n,k,4q+2}=u_{n,k,4q+3}
\end{equation}
for all $q<2^n-1$.
Set $C=(f|Q)^{-1}(Q_1)$,   $v_{n,k,p}=(h|C)^{-1}(u_{n,k,p})$, $$y_{n,k}=v_{n,k,0}^{-1}v_{n,k,1}v_{n,k,2}^{-1}\cdots v_{n,k,2^{n+2}-1},$$
and $Y=\{e_{F\,X}\}\cup\{y_{n,k}:n\in\w, k\geq k_n\}$.
Since $h|C:C\to\partial\J$ is a homeomorphism, the sequence $(y_{n,k})_{k\in\w}$
converges to  $e_{F\,X}$  for all $n\in\w$ by (\ref{convpre}).
In addition, the continuous homomorphism $Fh:F\,X\to F\,I$ maps $y_{n,k}$ to $r_{n,k}$ by our
choice of $v_{n,k,p}$, and hence $Fh(Y)=R$. By the definition, $S_\w$ is a union of a countable family of disjoint convergent sequences whose limit points coincide endowed with the strongest topology in which these
sequences are still convergent. 
Thus the continuity of $Fh|Y$ implies that $Y$ is homeomorphic to $R$.
\hfill $\Box$
\medskip

\noindent\textbf{Proof of Theorem~\ref{cor}.} 
Follows immediately from Proposition~\ref{main}, Lemma~\ref{tiur},
and the fact that $[0,1]$ is an absolute retract in the category of Tychonov spaces.
\hfill $\Box$
\medskip

\begin{lemma} \label{clcopy}
Under the assumptions of Corollary~\ref{coraa7} the image $i_X(X)$
is closed in $FX$.
\end{lemma}
\begin{proof}
Throughout the proof we shall identify $Z\in\Ob(\mathcal T)$
with $i_Z(Z)$. Let $j:X\to\bar{X}$ be the inclusion. It sufficies to show that
 $X=(Fj)^{-1}(\bar{X})$. Assuming the converse,
by the minimality of $F$ we could 
 find a finite subset $\{x_i:i\leq n\}\subset X$ and integers $m_i$, $i\leq n$,
such that $x^\ast:=(Fj)(x_0^{m_0}\cdot\cdots\cdot x_n^{m_n})\in \bar{X}\setminus X$. 
Therefore the elements  $x^\ast$ and $x_0^{m_0}\cdot\cdots\cdot x_n^{m_n}$
of $F\bar{X}$ coincide. Let $f:\bar{X}\to [0,1]$ be a continuous map such that
$f(x^\ast)\neq f(x_i)$ for all $i\leq n$.
From the above it follows that 
$$ f(x^\ast)=Ff(x^\ast)=Ff(x_0^{m_0}\cdots x_n^{m_n})=
f(x_0)^{m_0}\cdot\cdots\cdot f(x_n)^{m_n}, $$
which contradicts Lemma~\ref{injection}.
\end{proof}

\noindent\textbf{Proof of Corollary~\ref{coraa7}.}
Let  us denote by $\mathbb L$ the subspace 
$\{(0,0)\}\cup\{(\frac{1}{n},\frac{1}{mn}):n,m>0\}$ of $\IR^2$.
According to \cite[Lemma~8.3]{vD}, a first countable space contains a closed topological
copy of $\mathbb L$ if and only if it fails to be locally compact.

Assume, contrary to our claim, that $X$ is not scattered and fails to be locally compact.
Then $X$ contains a topological copy of the space $\IQ$ as well as a closed topological copy of
$\mathbb L$. Since $FX$ is generated by its second-countable subspace $X$,
it has countable pseudocharacter. In addition, it is normal being
Lindel\"of and Tychonov. 
Applying Theorem~\ref{cor}, we conclude that $FX$ contains a topological copy of
$S_\w$. It is well-known \cite{Na} that a topological group $G$ contains a (closed) topological copy of
$S_\w$ if and only if it contains a (closed) topological copy of the Arens' space $S_2$, see
\cite{Na} for its definition. It was shown by Lin \cite[Corollary~2.6]{Lin} that a regular space $Z$ with countable pseudocharacter contains a topological copy of $S_\w$
if and only if it contains a closed topological copy of $S_\w$.
It follows from the above  that $FX$ contains a closed topological copy of $S_\w$.
In addition, $FX$ contains a closed copy of $\mathbb L$ by Lemma~\ref{clcopy}, which
 contradicts \cite[Theorem~1]{Ba} 
asserting that a normal $k$-space containing closed topological copies of $\mathbb L$
and $S_\w$ is not homeomorphic to any topological group.
\hfill $\Box$
\medskip

Taras Banakh, Instytut Matematyki, Akademia \'Swi\c etokrzyska, Kielce (Poland) and \\
Department of Mechanics and Mathematics, 
 Ivan Franko Lviv National University,   Universytetska 1, Lviv, 79000, Ukraine.
\smallskip

\noindent \textit{E-mail address:}    \texttt{tbanakh@franko.lviv.ua}\\
\noindent \textit{URL:}   
\texttt{http://www.franko.lviv.ua/faculty/mechmat/Departments/Topology/index.html}

\medskip

 Du\v{s}an Repov\v{s}, Institute of  Mathematics, Physics and
Mechanics, Jadranska 19, P.O.B. 2964, Ljubljana, Slovenija 1001.
\smallskip

\noindent \textit{E-mail address:}   \texttt{dusan.repovs@guest.arnes.si}\\
\noindent \textit{URL:}   \texttt{http://pef.pef.uni-lj.si/$\tilde{\ }$dusanr/index.htm}

\medskip

 Lyubomyr Zdomskyy, Department of Mechanics and Mathematics,  Ivan Franko Lviv National University, Universytetska 1, Lviv, 79000, Ukraine.
\smallskip

\noindent \textit{E-mail address:}   \texttt{lzdomsky@gmail.com}
\end{document}